\documentstyle[12pt]{article}
\title{\bf On the Lie's Theorem for Lie  Color Algebras}
\author{Shouchuan Zhang $^{a,~b}$  \\
$a$. {\em Department  of Mathematics,
Hunan University } \\ {\em Changsha  410082, \
 P.R. China} \\
$b$. {\em Department of Mathematics, University of Queensland } \\
{\em Brisbane 4072, Australia } \\
 }
\date{}
\begin{document}
\newtheorem{Theorem}{\quad Theorem}[section]
\newtheorem{Proposition}[Theorem]{\quad Proposition}
\newtheorem{Definition}[Theorem]{\quad Definition}
\newtheorem{Corollary}[Theorem]{\quad Corollary}
\newtheorem{Lemma}[Theorem]{\quad Lemma}
\newtheorem{Example}[Theorem]{\quad Example}
\maketitle \addtocounter{section}{-1}

 \begin {abstract} We show that the Lie's Theorem holds for Lie color algebras
with a torsion-free abelian group $G$. We give an example to
show that the torsion-free condition is  necessary.
\end {abstract}

\section {Introduction}

Since Kac's work \cite{Ka77} on the classification of Lie superalgebras,
Lie color algebras have attracted much attention of mathematicians and
mathematical physicists. For instance, Scheunert \cite {Sc79} studied
Lie color algebras and
obtained the PBW basis and Ado theorem of these algebras. Montgomery
et al \cite {BFM96}
\cite {BFM01} studied the properties of Lie color algebras
determined by braided algebras,  $H$-Lie algebras (or generalized
Lie algebras) and $(H, r)$-Lie algebras. Majid \cite
{Ma94b}\cite {Ma00} introduced the basic theory of Lie braided
algebras and studied  the transmutation and bosonisation of Lie
braided bialgebras.

As generalized algebraic structures, Lie color algebras have found many
physical applications, e.g. in the construction of the Yang-Baxter
equation and integrable models \cite{McA97}\cite{McA98}.
In this letter we obtain the Lie's Theorem
for Lie color algebras. We show that the Lie's Theorem for Lie color
algebras holds provided
that the abelian group $G$ is  torsion-free. We give an
example to show that the torsion-free condition is  necessary. Our
results should shed new light on the study of structures of Lie color
algebras.

\section {Main Results} \label {e4}

The concepts of nilpotent elements, ad-nilpotent elements,
nilpotent Lie color algebras  and solvable Lie color algebras are
similar to ones in Lie algebras.

Let $G$ be an abelian group with a skew symmetric bicharacter $r$.
If $V = \sum _{g\in G}\oplus V_g$ is a
graded vector space over field $F$ with $dim \ V_g = n_g < \infty$
and $A_{gh}=Hom_{F}(V_h, V_g)$ for any $g, h \in G$, then
$A=\sum\{A_{ij} \mid i, j\in G \}$ is an algebra graded by $G$ with
$A_g=\sum_{i=j+g}\,A_{ij}$ for any $g\in G$. $A$ is a Lie color algebra
under the bracket operation $[a,b]=ab-r(|b|,|a|)ba$ for any homogeneous
elements $a,b\in A$. Here $|a|$ stands for the degree of the element
$a$. This Lie color algebra is referred to as
general linear Lie color  algebra, written as $gl (\{n_g\}, F)$
or $gl (\{V_g\}, F)$. Its Lie color subalgebras are called
linear  Lie color algebras. In fact, $gl (\{n_g\}, F) = \{ f
\in End _k V \mid
 ker f \hbox { is finite codimensional } \}.$ When $G$ is finite,
 we may view $gl (\{n_g\}, F)$ as a
block matrix algebra over $F.$ When $G=0$, we denote $gl (\{n_g\},
F)$ by $gl (n, F)$, in this case it is the ungraded  general
linear Lie algebra.


\noindent
\begin {Lemma}\label {4.1} Let  $L$ be a Lie color algebra and $N$
be a Lie color subalgebra of $gl(\{V_g\},F)$. Then,
(i)  $(ad \ X)^m (Y) = \sum _{i +j = m} k_{ij} X^i Y X^j$ for any
homogeneous elements  $X, Y \in  N,  k_{ij } \in F,$
 and natural number $m$.

(ii) If  $X \in N$  is nilpotent homogeneous, then $ad \ X$ is
nilpotent too.

(iii) For any homogeneous  $x, y, z \in L$,  $[[x, y], z] = [x [y,
z]] + r (\mid z \mid , \mid y \mid )[ [x, z ] , y] $, i.e. $ad\
[x, y] = [ad \ x, ad \ y].$ In other word, $ad $ is a color
representation of Lie color algebra $L.$

(iv) $Z(L)$ is a color  ideal of Lie color algebra $L$.
Here $Z(L) = \{ x \in L \mid  [x , y] = 0  \hbox { for any
} y \in L \}$, called the center of $L$.
\end {Lemma}

\noindent {\bf Proof.} (i) It can be obtained   by induction on natural
number $m$.

(ii)    It immediately follows from (i).

(iii) It is trivial.

(iv) We only need show that $Z(L) $ is $G$-graded. For any $x \in
Z(L)$ and for any homogeneous element $ y  \in L$,  since $[x ,
y]=0$, we have that $\sum _{g\in G}[ x_g, y ] =0 ,$ which implies
that $[x_g, y]=0$. Consequently, $x_g \in Z(L).$  \  $\Box$

\noindent
\begin {Theorem}\label {4.2} Let $L$ be a finite-dimensional
Lie color subalgebra of $gl (\{V_g\}, F).$ If every homogeneous
element in $L$ is nilpotent and $V \not=0,$ then there exists  a
non-zero homogeneous element $v\in V$ such that $Lv=0.$
\end {Theorem}

\noindent {\bf Proof.}  We use induction on $dim \ L = n.$  The case $dim \
L = 0 $ is  obvious.  Let $K$ be a color subalgebra of $L$ with $K
\not= L$. For any $X \in K$, it is clear that there exists an
$F$-linear map $\phi (X)$ from $L/K$ to $L/K$ such that the below
diagram becomes  commutative.
$$\begin {array}{lcr}
{}&            ad \ X &{} \\
 L & \longrightarrow&  L   \\
\pi   \downarrow & \searrow & \downarrow    \pi \\
 L/K  & \longrightarrow  & L /K  \\
{}&           \phi (X) &{} \\
\end {array},  $$ i.e.  for any $Y\in L,$ $\phi (X) (Y+K) = [X,Y] +K. $
It is easy to verify that $\phi $ is a Lie color algebra
homomorphism from $K$ to $gl (\{L_g/K_g\},F) $.
 The induction  hypothesis guarantees
existence of a homogeneous element  $X_0 +K \not= K$ in $L/K$ such
that $\phi (X)(X_0+K) = K $, i.e. $[X, X_0] + K = K$ and $[X, X_0]
\subseteq K $ for any $X \in K.$ Thus $K + F X_0$ is a Lie color
subalgebra of $L$ with $dim (K+ FX_0)= dim K +1.$ Consequently,
there exist a Lie color subalgebra $K$ and a homogeneous element
$z\in L$ such that $L = K + F z.$ Obviously,  $K$ is a color ideal
of $L.$

By induction, $W = \{ v \in V \mid Kv =0   \}\not= 0$. Now we show
that $W$ is a graded space. For any homogeneous element  $0\not= x
\in K$, $0\not= w\in W$, we see that $x(w) = \sum _{g\in G} x(w_g)
=0$ implies $x(w_g) = 0$ since $x$ is homogeneous. Thus $w_g \in
W$ for any $g\in G,$ which implies that $W$ is a graded space.
Since $K$ is a $G$-graded ideal,  $W$ is stable under $L$: for
three homogeneous elements $x\in L , y \in K, w\in W$, we see that
$y(x(w)) =[y, x] (w) + r(\mid x\mid , \mid y \mid ) x(y (w))=0. $
Since $z$ is nilpotent on $W$ there exists a  non-zero $v \in W$
such that $z(w)=0$ and $\sum _{g\in G} z(w_g)=0$, which implies
$z(w_g)=0$ for any $g\in G.$ Consequently, there exists a non-zero
homogeneous $w_g \in W$ such that $z(w_g)= 0$. Finally, $L(w_g) =
0 $, as desired. $\Box$ In \cite{Sc79a}, Scheunert stated that the
Engel's theorem is true by the usual proof. Here for completeness,
we have written down the theorem and provided a short proof.
\noindent
\begin {Theorem}\label {4.3}  (See \cite {Sc79a})(Engel's Theorem)
Let $L$ be a finite-dimensional Lie color algebra .
If all homogeneous elements of $L$ is $ad$-nilpotent , then $L$ is
nilpotent.
\end {Theorem}

\noindent {\bf Proof.} It is clear that $ad \ L =\{ad \ x \mid  x \in L\}
\subseteq gl (\{L_g \}, F) $ satisfies the hypothesis of Theorem
\ref {4.2}. Thus there exists  a homogeneous element
$0\not= x $ in $L$ for which $[L,x] =0, $ i.e. $Z(L)\not=0$. Now
$L/Z(L)$ evidently consists of $ad$-nilpotent elements and has
smaller dimension than $L$. Using induction on $dim \ L$, we find
that $L/Z(L) $ is nilpotent, which implies $L$ is nilpotent.
$\Box$

\noindent
\begin {Lemma}\label {4.7} Let $A\in gl(\{n_g\}, F)_u$ with
torsion-free $u\not=0$. Then

(i)  If $A$ has an eigenvalue,  then  $A$ must have a zero
eigenvalue and has homogeneous eigenvector of eigenvalue zero.

(ii) every  eigenvalue  of $A$ is zero.

(iii) $A$ is nilpotent.
\end {Lemma}

\noindent {\bf Proof.} (i) If $\lambda $ is an eigenvalue of $A$, then there
exists an eigenvector $v \in V$ such that $A(v) = \lambda v$, i.e.
$$\sum _{g\in G} Av_g = \sum _{g\in G} \lambda v_g.$$ Thus $A(v_g) =
\lambda v _{g + u}$ for any $g\in G.$  Assume $v_{g_0} \not=0.$
There exists a natural number $m$ such that $v_{g_0 + mu} \not=0 $
and $v_{g_0 + (m+1) u}=0.$ Since $A(v_{g_0 + m u}) =\lambda v_{g_0
+ (m+1) u}=0$, which implies that $A$ has a zero eigenvalue and
this eigenvector is homogeneous.

(ii) If $\lambda $ is a non-zero eigenvalue of $A$, then there
exists an eigenvector $v \in V$ such that $A(v) = \lambda v$, i.e.
$$\sum _{g\in G} Av_g = \sum _{g\in G} \lambda v_g.$$ Thus $A(v_g) =
\lambda v _{g + u}$ for any $g\in G.$  Assume $v_{g_0} \not=0.$
There exists a natural number $m$ such that $V_{g_0 - mu} \not=0 $
but $v_{g_0 - (m+1) u}=0.$ Since $A(v_{g_0 - (m+1) u}) =\lambda
v_{g_0 - m u}=0$,  $\lambda =0.$ We get a contradiction.

(iii) follows from part (ii). $\Box$

\noindent
\begin {Theorem} \label {4.8}
Let $L$ be a finite-dimensional solvable Lie color
subalgebra of $gl ( \{V_g \}), F)$ over algebraically closed field
$F$ and all of homogeneous elements in $[L,L]$ be nilpotent. If
$V\not=0$ and $G$ is a torsion-free abelian group , then $V$
contains a common homogeneous eigenvector for all the
endomorphisms in $L.$
\end {Theorem}
{\bf Proof.} Use induction on $dim \ L.$ The case $dim \ L =0$ is
trivial. We attempt to imitate the proof of Theorem \ref {4.2}.
The idea is  (1)  to locate a color ideal $K$ of codimension one,
(2) to show by induction that common homogeneous eigenvectors
exist for $K$, (3) to verify that $L$ stabilizes a space
consisting of such eigenvectors , and (4) to find in that space
homogeneous eigenvector for a single homogeneous $x\in L$
satisfying $L = K + Fz.$

Step (1)    is easy. Since $L$ is solvable, $[L, L] \not= 0$ and
$L$ properly includes $[L, L]$.  $L/ [L, L] $  being abelian, any
$G$-graded subspace is automatically an ideal. Take a $G$-graded
subspace  of codimension one, then its inverse image $K$ is a
color ideal of codimension one in $L$. Set $ L = K + Fz$ for some
homogeneous element $z \in L.$

 For step (2), use induction to find a common homogeneous
 eigenvector $v \in V$
for  $K$. This means that for $x\in K$, $x (v) = \lambda (x) v,$
$\lambda : K \rightarrow F $ is some linear function. Fix the
$\lambda$ and denote by $W$ the color subspace generated by
$$ \{w \in V \mid  w  \hbox { is homogeneous and } x (w) = \lambda (x) w,
\hbox  { for all } x \in K \};
 \hbox { so } W \not= 0.$$

 Step (3) consists in showing that $L$ leaves $W$ invariant.
 Assuming for the moment that this is done, preceed to step (4);
 write $L= K +Fz$  and use the fact  that $F$ is algebraically
 closed and Lemma \ref {4.7}
 to find a homogeneous  eigenvector $v_0\in W$ (for some eigenvector of $z$).
Then $v_0$ is obviously a common eigenvector for $L$ (and $\lambda
$ can
 be extended to a linear function on $L$ such that
 $x(v_0) = \lambda (x) v_0, x\in L$).

 It remains to show that $L$ stabilizes $W$. Let  homogeneous $w\in
 W$, homogeneous $x\in L$. To test whether or not $x (w)$ lies in  $
 W$, we must take arbitrary  $y \in K$ and examine
$y (x(w)) = [y,x] (w) + r(\mid x \mid , \mid y \mid)x(y(w)) =
\lambda ([y,x] )(w) + r(\mid x \mid , \mid y \mid) \lambda
(y)x(w).$ Thus we have to  prove that $\lambda ([y,x ])=0$.
 However, since $[y,x ]$ nilpotent, $\lambda ([y,x ])=0$. We complete the
 proof. $\Box$

\noindent
\begin {Corollary} \label {4.9}  (Lie's Theorem )
Let $L$ be a finite-dimensional
solvable Lie color subalgebra of $gl (\{V_g\}, F)$ and all of
homogeneous elements in $[L,L]$ nilpotent. If $G$ is a
torsion-free abelian group, then $L$ stabilizes some color flag in
$V$(in other words, the matrices of $L$ relative to a suitable
homogeneous  basis of $V$ are uppertriangular )  .
 \end {Corollary}

The Lie's theorem is one of the main results in the present letter.
The condition  that   $G$ is  torsion-free is necessary  for the Lie's
theorem to be true as is seen from the following example.

\begin {Example} \label {4.10}
Let $G = {\bf  Z}_3$ with  trivial bicharacter $r(x,y) = 1$ for $x, y \in {\bf  Z}_3$,  $\{e_1, e_2, e_3\}$  the
natural basis of ${\bf F}^3$,  $V_i = Fe_i$ for $i = 1, 2, 3$ and
a linear transformation $   A = \left (  \begin{array}{l l l}0 & 1 & 0\\
0 & 0 & 1\\
1 & 0 & 0\\
\end{array} \right )$.
Let $L = F A$ be a linear Lie color algebra in $gl (\{V_g\}, F)$.
 Obviously, the degree of $A$ is $2$ and $[L,L]=0$.
In  other word,   $L$ is  a finite-dimensional solvable color
subalgebra of $gl (\{V_g\},F)$ and all  homogeneous elements in $[L,L]$
are nilpotent. But  the matrices of linear transformation $A$
relative to any homogeneous  basis of $V$ can not be upper triangular.
 \end {Example}

\noindent {\bf  Remark:} In the case of Lie algebras,
the other conditions in Corollary  \ref {4.9} imply the
condition that every homogeneous element in  $[L,L]$ is nilpotent by
 \cite [Corollary 4.1 A] {Hu72}.
In other word, our Corollary \ref {4.9} implies the Lie's
theorem for (ungraded) Lie algebras.

\noindent
\begin {Corollary} \label {4.11}    Let $L$ be a finite-dimensional solvable
 Lie color subalgebra of $gl (\{ n_g\}, F)$ and all  homogeneous elements in
 $[L,L]$be
nilpotent. If  $G$ is a torsion-free abelian group,  then there
exists a chain of color ideals of $L$, $0 = L_0 \subseteq L_1
\subseteq \cdots \subseteq L_n = L$ such that $dim L_i = i.$
 \end {Corollary}

\noindent
{\bf Acknowledgement }: The work is supported by Australian Research Council.
S.C.Z thanks the Departemnt of Mathematics, University of
Queensland for hospitality.

\end {document}